\documentclass[12pt]{amsart}
\usepackage{amsmath,amssymb}
\numberwithin{equation}{section}
\theoremstyle{plain}
  \newtheorem{theorem}{Theorem}[section]
  
  \newtheorem{proposition}[theorem]{Proposition}
  \newtheorem{lemma}[theorem]{Lemma}
  \newtheorem {corollary}[theorem]{Corollary}
 
\theoremstyle{definition}

\theoremstyle{remark}
 \newtheorem{remark}[theorem]{Remark}

\newcommand{\op}[1]{{{}^\sigma\!{#1}}}
\newcommand{\In}[1]{\daffH\*\nolimits_{\ratH}{#1}}
\def\H{\bf H}

\def\affWreg{{\affW^\circ}}

\def\nil{{\mathfrak N}}

\def\rat{{{}}}

\def\af{{\rm aff}}
\def\daffH{{\widetilde\H_\kappa}}
\def\daffHm{{\widetilde\H_{-\kappa}}}
\def\daffL{\widetilde{\mathcal L}_\kappa}
\def\daffst{\widetilde{\Delta}_\kappa}
\def\daffHmod{{\daffH{\hbox{-}}{\mathcal Mod} }}
\def\ratL{{\mathcal L}^\rat_\kappa}
\def\ratst{{\Delta}^\rat_\kappa}
\def\ratHmod{{ \ratH{\hbox{-}}{\mathcal Mod} }}
\def\ratH{{\H_\kappa}}
\def\ratHp{{\H_{\geq0}}}

\def\affH{{\H^\af}}

\def\Preg{P_\circ}

\def\O{{\mathcal O}}

\def\Fr{{\mathcal F}}
\def\xch{y}
\def\Xch{\F[\underline{y}]}
\def\X{\F[\underline{x}]}
\def\XL{\F[\underline{x}^{\pm1}]}
\def\ev{\psi}
\def\emb{\iota}
\def\embinv{\jmath}
\def\isomto{\stackrel{\sim}{\to}}

\def\part{{\Lambda}^+}

\def\affWreg{{\aff {\mathfrak S}_n^\circ}}

\def\Oss{{{\mathcal O}^{ss}}}

\def\Q{{\mathbb Q}}
\def\Z{{\mathbb Z}}
\def\F{{\mathbb F}}

\def\W{{\mathfrak S}_n}
\def\affW{{\aff{\mathfrak{S}}_n}}

\def\e{{\epsilon}}
\def\ech{u}

\def\ep{\epsilon}
\def\ee{x}
\def\UU{\F[\underline{u}]}
\def\lm{{\lambda}}

\def\m{{m}}
\def\t{{t}}
\def\wt{{\F^n}}
\def\Res{{{\rm{Res}}_{\ratH}}}

\def\End{{{\rm{End}}}\,}
\def\Hom{{{\rm{Hom}}}\,}

\def\Irr{{\hbox{\rm{Irr}}}}
\def\dim{{\hbox{\rm{dim}}}}
\def\id{{\hbox{\rm{id}}}\,}

\def\Coker{{\rm Coker}}

\def\aff{\widehat}

\def\qed{\hfill$\square$}

\def\+{\mathop{\oplus}}
\def\*{\mathop{\otimes}}
\def\induction{\daffH\*_{\ratH}(-)}

\newcommand{\omitted}[1]{}
\begin{document}
\title[Double affine Hecke algebras]{Rational and
trigonometric degeneration
of the double affine Hecke algebra of type $A$}
\author{Takeshi Suzuki}
\address{Research Institute for Mathematical Sciences\\
Kyoto University, Kyoto, 606-8502, Japan}
\email{takeshi@kurims.kyoto-u.ac.jp}
\begin{abstract}
We study a connection between
the representation theory of
the rational Cherednik algebra of type $GL_n$ 
and
the representation theory of
the degenerate double affine Hecke algebra (the degenerate DAHA).
We focus on an algebra embedding from the rational Cherednik algebra
to the degenerate DAHA 
and investigate the induction functor through this embedding.
We prove that this functor
embeds the category ${\mathcal O}$ for the rational Cherednik algebra
fully faithfully
into the category ${\mathcal O}$ for the degenerate DAHA.
We also study the full subcategory ${\mathcal O}^{ss}$ of ${\mathcal O}$
consisting of those modules which
are semisimple with respect to the commutative subalgebra
generated by Cherednik-Dunkl operators.
A classification of all irreducible modules in ${\mathcal O}^{ss}$
for the rational Cherednik algebra
is obtained from the corresponding
result for the degenerate DAHA.
\end{abstract}
\maketitle
\begin{center}
  {\bf Introduction}
\end{center}
We study a connection between
the representation theories of 
two algebras of type $GL_n$;
the rational Cherednik algebra 
$\ratH$ and
the degenerate double affine Hecke algebra
(or the trigonometric Cherednik algebra) 
$\daffH$.
These two algebras can be thought (see \cite{EG}) as a rational
and trigonometric degeneration of
the double affine Hecke algebra 
introduced by Cherednik
\cite{Ch;unification, Ch;double} respectively.

The categories of finite dimensional representations for
the rational and trigonometric Cherednik algebras 
(of type $SL_n$) have been related
through deformation argument 
(\cite[Proposition 7.1]{BEG}\cite{Ch;fourier}).

In this paper, we relate what is called the categories $\O$ 
(\cite{DO,GGOR,VV}) for
$\ratH$ and $\daffH$ by a different approach.
We give an algebra embedding from
 $\ratH$ to $\daffH$,
which extends to an algebra isomorphism from
 a localization of $\ratH$ to $\daffH$ 
(Proposition~\ref{pr;embedding}),
and investigate the induction functor
$\daffH\*_{\ratH}(-):\ratHmod\to\daffHmod$
through this embedding.

By restricting this functor, we obtain the functor
from $\O(\ratH)$ into $\O(\daffH)$, where
 $\O(\ratH)$ and $\O(\daffH)$ denote the category $\O$ for
$\ratH$ and $\daffH$ respectively.
The main result of this paper (Theorem~\ref{th;fullyfaithful}) 
asserts that this functor
embeds $\O(\ratH)$ into $\O(\daffH)$ fully faithfully.
Moreover it follows
that our functor sends a standard module of $\ratH$ to 
a standard module of
$\daffH$, and an irreducible module to an irreducible module.

Through this functor, one can transform
results on the one side to the other.
For example, the multiplicity
of an irreducible module in the composition series
of a standard module for $\ratH$ is equal to
the corresponding multiplicity for $\daffH$.
Note that $\O(\daffH)$ has been also related to
the category $\O$ for the (original) double affine Hecke algebra
(\cite{VV}),
and hence this multiplicity is also equal to
the  corresponding multiplicity for the double affine Hecke algebra,
which can be expressed in terms of Kazhdan-Lusztig polynomials
as shown in \cite{Va} by a geometric method.

We also focus on the full subcategory $\Oss$ of $\O$
consisting of those modules which
are semisimple with respect to the commutative subalgebra
which is realized as the algebra generated by
Cherednik-Dunkl operators on the polynomial representation.
This class of representations for $\daffH$ have been
studied in \cite{Ch;fourier,SV} (see also \cite{Ch;introduction}),
and in particular the classification
of all irreducible modules in $\Oss(\daffH)$
has been obtained.
It turns out that the induction functor transforms
$\Oss(\ratH)$ into $\Oss(\daffH)$, and hence
the classification (Theorem~\ref{th;classification})
of the irreducible modules in $\Oss(\ratH)$
follows from the classification result for $\Oss(\daffH)$.

\medskip
\noindent
{\bf Acknowledgments.}
The author would like to thank Prof.\,Opdam
for discussion and his lecture on Cherednik algebras while
his stay in RIMS.
The author is also very grateful to Prof.\,Kashiwara
for important suggestion and comments.
%
\section{The affine Weyl group}
%
Let $\F$ be
an algebraically closed
field of characteristic $0$. 
Let $\Q$ and $\Z$ denote the set of rational numbers
and the set of integers respectively, which are subrings of $\F$.

Throughout this paper,
 we use the notation
$$[a,b]=\{a,a+1,\dots,b\}\subset \F$$
for $a,b\in\F$ such that $b-a\in\Z$.

Let $n\in\Z_{\geq 2}.$
Let $\W$
denote
the symmetric group of degree $n$.
For $i,j\in[1,n]$ with $i\neq j$,
let $s_{ij}$ denote
the transposition of $i$ and $j$.
In particular, we write  $s_{i}\ (i\in[1,n-1])$ for
the simple reflection which transposes $i$ and $i+1$.

Let $P=\oplus_{i\in[1,n]}\Z\e_i$
denote the free $\Z$-module generated by $\ep_1,\ep_2, \dots,\ep_n$,
on which $\W$ acts naturally.
Define {\it the extended affine Weyl group} $\affW$
as the semidirect product
$P\rtimes \W$.

Put $
\pi=\t_{\ep_1} s_{1}s_2\cdots s_{n-1}$
and $s_0=\t_{\ep_1-\ep_n}s_{1n}$,
where $\t_\eta$ denotes the element of $\affW$ corresponding to
$\eta\in P$.
Then $\{\pi, s_0,s_1,\dots,s_{n-1}\}$
is a set of generators
 of the group $\affW$. 
We put
$\Preg=\oplus_{i\in[1,n]}\Z_{\geq 0}\e_i\subset P$
and $\affWreg=\Preg\rtimes\W\subset \affW$.
Observe that $\affWreg$ is the semigroup (with unit) generated by
 $\{\pi, s_1,s_2,\dots,s_{n-1}\}$.

 For a semigroup $G$, let $\F G$ denote the semigroup algebra of $G$.
 \begin{lemma}\label{lem;pi_and_W}
   The linear map
$\F[\pi,\pi^{-1}]\*\F\affWreg \to \F\affW$
given by $\pi^k\* w\mapsto \pi^kw$
induces a linear isomorphism
$\F[\pi,\pi^{-1}]\*_{\F[\pi]}\F\affWreg \isomto \F\affW$.
 \end{lemma}
\noindent{\it Proof.} 
Let $w\in\affW$.
Note that $\pi^n=t_{\e_1+\e_2+\cdots+\e_n}$. 
Hence there exist
some $k\in\Z$ and $\bar w\in\affWreg$
such that $w=\pi^k \bar w$.
It follows that
the correspondence $w\mapsto \pi^k\*\bar w$
gives rise to a well-defined linear map
$\F\affW\to\F[\pi,\pi^{-1}]\*_{\F[\pi]} \F\affWreg$,
and it gives the inverse of the induced map
$\F[\pi,\pi^{-1}]\*_{\F[\pi]}\F\affWreg \to \F\affW$.
\qed

\medskip\noindent
We denote the polynomial ring and the Laurent polynomial ring
in $n$-variables $x_1,x_2,\dots,x_n$ by $\X$
and $\XL$ respectively:
\begin{align*}
  \X&=\F[x_1,x_2,\dots,x_n],\quad
\XL=\F[x_1^{\pm1},x_2^{\pm1},\dots,x_n^{\pm1}].
\end{align*}

In the sequel, 
we identify
$\F P$ with $\XL$, and identify $\F \affW$ with $\XL\rtimes \W$
as algebras
through the correspondence $t_{\e_i}\mapsto x_i$ $(i\in[1,n])$.
Observe that $\F\affWreg$ is identified with $\X\rtimes \W$.
\section{Cherednik algebras}
Fix $\kappa\in \F\setminus\{0\}$.

Let $\Fr^\rat$ denote the associative algebra
generated freely by the algebras $\F\affWreg=\X\rtimes \W$ and 
$\Xch=\F[y_1,y_2,\dots,y_n]$
such that the natural inclusions
$\F\affWreg\to\Fr^\rat$
and $\Xch\to\Fr^\rat$ are algebra homomorphisms.
{\it The rational Cherednik algebra} $\ratH$ of type $GL_n$
is the unital associative $\F$-algebra defined as
the quotient algebra of $\Fr^\rat$
by the following relations:
\begin{eqnarray*}
\label{eq;ratrel}
& &s_i y_j =y_{s_i(j)} s_i
\quad (i\in[1,n-1],\ j\in[1,n]),
\\
\label{eq;ratrel3}
& &\left[ \xch_i, x_j\right] =
\begin{cases}
\kappa+\sum_{i\neq k}s_{ik}\quad &(i=j),\\
-s_{ij}\quad &(i\neq j).
\end{cases}
\end{eqnarray*}
Let $\widetilde\Fr$ denote  the algebra generated freely by
the algebras $\F\affW=\XL\rtimes \W$ and $\UU=\F[u_1,u_2,\dots,u_n]$
such that the natural inclusions
$\F\affW \to\widetilde\Fr$
and $\UU\to\widetilde\Fr$ are algebra homomorphisms.
{\it The degenerate double affine Hecke algebra}
(or {\it the trigonometric Cherednik algebra}  )
$\daffH$ of type $GL_n$
is the unital associative $\F$-algebra defined as
the quotient algebra of $\widetilde\Fr$
by the following relations:
\begin{align*}
&s_i u_i =u_{i+1} s_i-1,\
s_i u_{i+1} =u_{i} s_i+1
\quad (i\in[1,n-1]),\\
&s_i u_j =u_j s_i
\quad (i\in[1,n-1],\ j\in[1,n],\ j\neq i,i+1),\\
&\left[ u_i, x_j\right]=
\begin{cases}
\kappa x_i+\sum_{1\leq k<i} {x_k}s_{ki}
+\sum_{i<k\leq n} {x_i}s_{ik}
\quad &(i=j),\\
-{x_j}s_{ji}\quad &(i>j),\\
-{x_i}s_{ij}\quad &(i<j).
\end{cases}
\end{align*}
\begin{remark}
The degenerate double affine Hecke algebra can be defined
in an alternative way
in terms of the generators
$$\{\pi,\pi^{-1},
s_0,\dots,s_{n-1},\ech_1,\dots,\ech_n\}.$$
It is well-known
that the two definitions coincide
(see e.g. \cite[Proposition 1.3.7]{AST}).
Note in particular that the following relations hold:
\begin{align}\label{eq;pi_ech}
\pi \ech_i&=\ech_{i+1}\pi\ (i\in[1,n-1]),\quad
\pi \ech_n=(\ech_1-\kappa)\pi.
\end{align}
\end{remark}
%

\begin{proposition}
\label{pr;PBW}
$({\rm Cherednik}$, \cite[Theorem 1.3]{EG}$)$

\smallskip\noindent
${\rm (i)}$ Multiplication in $\ratH$ induces a
linear isomorphism
$$\X\* \F\W \* \Xch\overset{\sim}{\rightarrow}\ratH.$$

\smallskip
\noindent
${\rm (ii)}$  Multiplication in $\daffH$ induces a
linear isomorphism
$$\XL\* \F\W \* \UU\overset{\sim}{\rightarrow}
\daffH.$$
\end{proposition}
The subalgebra $\affH=\F\W\cdot \UU$
of $\daffH$ is called {\it the degenerate affine Hecke algebra}
of type $GL_n$.
%
\section{The algebra embedding and the induction functor}
The algebra structure of $\ratH$ is extended to the
algebra structure of $\XL \*_{\X} \ratH=\XL\*\F \W\* \Xch$.
%
\begin{proposition}\label{pr;embedding}
${\rm (i)}$ There exists a unique algebra homomorphism
$\emb : \ratH \to \daffH$
such that
\begin{align*}
\emb(w)&=w \quad (w\in\W),\\
\emb(x_i)&=\ee_i \quad (i\in[1,n]),\\
\emb(\xch_i)&=x_i^{-1}\left(\ech_i-\sum\nolimits_{1\leq
j<i}s_{ji}\right)\quad
(i\in[1,n]),
\end{align*}
and moreover it is injective.

\smallskip
\noindent
${\rm (ii)}$
The algebra embedding $\emb:\ratH\to\daffH$
is extended to the algebra isomorphism
 $\XL\*_{\X} \ratH\isomto
\daffH$, whose inverse map
$\embinv:\daffH\to\XL\*_{\X} \ratH$
is defined by 
\begin{align*}
\embinv(w)&=w \quad (w\in\W),\\
\embinv(\ee_i)&=x_i \quad (i\in[1,n]),\\
\embinv(\ech_i)&=
x_i\xch_i+\sum\nolimits_{1\leq j<i}s_{ji}\quad(i\in[1,n]).
\end{align*}
\end{proposition}
\noindent
{\it Proof.}
By checking the defining relations of $\ratH$ (resp. $\daffH$)
by  direct calculations,
one can show that the assignment
$\emb$ (resp. $\jmath$) are uniquely
extended to the algebra homomorphism $\XL\*_{\X} \ratH\to\daffH$
(resp. $\daffH\to\XL\*_{\X} \ratH$).
It is easy to see that $\embinv\circ\emb=\id_{\daffH}$
and
$\emb\circ\embinv=\id_{\XL\*_{\X} \ratH}$.
This proves (ii).
Now (i) follows immediately
because the natural map $\ratH
\to\XL\*_{\X} \ratH$
is injective by Proposition~\ref{pr;PBW}-(i).
\qed
\begin{remark}
Recall the polynomial representation
$\ratH\to\End_\F(\X)$, where
the Dunkl operators
$$Y_i=\kappa\frac{\partial}{\partial x_i}+\sum_{j\neq i}
\frac{1}{x_i-x_j}(1-s_{ij})\ \ (i\in[1,n])$$
appear as the images of $\xch_i$ $(i\in[1,n])$.
The operators 
$U_i:=
x_i Y_i+\sum\nolimits_{1\leq j<i}s_{ji}$
$(i\in[1,n])$, which are also pairwise commutative, 
are introduced by Cherednik (\cite{Ch;unification})
 as a trigonometric 
analogue of
the Dunkl operators. 
The original and trigonometric Dunkl operators
play important role in the theory of orthogonal polynomials
(see e.g. \cite{DX}).
Note in particular that
the simultaneous eigenvectors of $\X$
with respect to the trigonometric Dunkl operators
 $U_1,\dots,U_n$ are
the nonsymmetric Jack polynomials.

Note also that generalization of trigonometric Dunkl operators
to other types than $A$ has been
 studied (\cite{K,DO}).
\end{remark}
\begin{remark}
The restriction of the isomorphism $\embinv$ to the subalgebra $\affH$
gives an algebra embedding $\affH\to\ratH.$
More generally, it was shown in \cite[Proposition 4.3]{EG}
that
the assignment
$$\ech_i\mapsto ax_i+b\xch_i+x_i\xch_i+\sum_{1\leq j<i}s_{ji}
\ \ (i\in[1,n]),\ \
w\mapsto w\ \ (w\in\W)$$
extends to an algebra embedding $\affH\to\ratH$
for any $a,b\in \F$.
When $b=0$, this embedding extends to the
algebra isomorphism $\embinv_a:\daffH\to \XL\*_{\X}\ratH$
such that $\embinv_a( \ee_i)=x_i$ $(i\in[1,n])$.
\end{remark}
In the rest of the paper,
we often identify $\daffH$ with $\XL\otimes_{\X}\ratH$,
and regard $\ratH$ as a subalgebra of $\daffH$.

Let us consider the induction functor
$M\mapsto \daffH\*_\ratH M$
from the category $\ratHmod$ of $\ratH$-modules to
the category $\daffHmod$ of $\daffH$-modules.
Since $\daffH$ can be seen as a
localization of $\ratH$
by Proposition~\ref{pr;embedding},
we have
\begin{corollary}\label{cor;flat}
The $\daffH$ is flat as an $\ratH$-module,
or equivalently, the induction functor
 $\daffH\otimes_{\ratH}(-):\ratHmod\to\daffHmod$
is exact.
\end{corollary}
\begin{corollary}\label{cor;ind_irr}
Let $M$ be an irreducible $\ratH$-module.
Then $\daffH\*_\ratH M$ is an irreducible $\daffH$-module
or zero.
\end{corollary}
\noindent{\it Proof.}
Follows easily from the 
isomorphism
 $\daffH\*_\ratH M\cong \XL\*_{\X} M$.
\qed

\medskip
For later use,
we give another expression of $\daffH$, which follows from 
Lemma~\ref{lem;pi_and_W}.
%
\begin{corollary}\label{cor;pi_and_H}
The linear map $\F[\pi,\pi^{-1}]\* \ratH\to \daffH$
given by $\pi^k\* v\mapsto \pi^kv$.
induces a linear isomorphism
$\F[\pi,\pi^{-1}]\*_{\F[\pi]} \ratH\isomto \daffH.$
\end{corollary}
We will use the following lemma later.
\begin{lemma}
There exists an algebra isomorphism
$\sigma:\daffH\to \daffHm$
such that
$$\sigma(s_i)=-s_i\ (i\in [0,n-1]),\
\sigma(x_i)=(-1)^{n-1}x_i,\ 
\sigma(\ech_i)=-\ech_i\ (i\in[1,n]),
$$
and its restriction to $\ratH$ gives an isomorphism
$\sigma: \ratH\isomto \H_{-\kappa}$.
\end{lemma}
%
\section{The category $\O$}
Put $\Xch_+=\sum_{i\in[1,n]}\xch_i\Xch\subset \ratH$.
An element $v$ of an $\ratH$-module $M$
is said to be $\Xch_+$-nilpotent
if there exists $k>0$ such that
$\xch_i^k v=0$ for all $i\in[1,n]$.
An $\ratH$-module $M$ is said to be
locally nilpotent for $\Xch_+$, if
any element of $M$ is  $\Xch_+$-nilpotent.
Let $\O(\ratH)$ denote the full subcategory of $\ratHmod$ consisting of
those modules which are finitely generated over $\ratH$ and
locally nilpotent for $\Xch_+$.

Let $\O(\daffH)$ denote the  full subcategory of $\daffHmod$
consisting of those modules which are
finitely generated over $\daffH$ and locally finite for $\UU$, i.e.,
a finitely generated $\daffH$-module $M$ is in $\O(\daffH)$
if $\UU v<\infty$ for any $v\in M$.
\begin{proposition}\label{pr;finitecomp} 
$\rm{(i)}$
$($\cite[Corollary 2.16]{GGOR}$)$
The category $\O(\ratH)$ is  a Serre subcategory of $\ratHmod$,
and
any object in $\O(\ratH)$ has a finite composition series.

\smallskip\noindent
$\rm{(ii)}$ $($\cite[Proposition 2.1, Proposition 2.2]{VV}$)$
The category $\O(\daffH)$ is  a Serre subcategory of $\daffHmod$,
and
any object in $\O(\daffH)$ has a finite composition series.
\end{proposition}
For $d\in\Z_{\geq0}$, define $\X_{(d)}$
as the subspace of $\X$
spanned by $\{x_1^{d_1}x_2^{d_2}\cdots x_n^{d_n}\mid
\sum_{i\in[1,n]}d_i=d\}$.
We have $\oplus_{d=0}^\infty\X_{(d)}=\X$.
Define $\Xch_{(d)}$ similarly.
\begin{proposition}\label{pr;OtoO}
Let $M$ be an object of $\O(\ratH)$.
Then $\daffH\*_\ratH M$ is an object of $\O(\daffH)$.
\end{proposition}
\noindent{\it Proof.}
Recall that, under the identification $\daffH\cong\XL\*_{\X}\ratH$,
the element $u_i\in\daffH$ is expressed as $x_iy_i+\sum_{j<i}s_{ji}$,
and it is contained in $\ratH\subset \daffH$ .

Let $M$ be an object of $\O(\ratH)$ and take $v\in M$.
As a step, we will prove that
$\UU v\subseteq M$ is finite-dimensional.
It is easy to see that
$$\ech_i^d=\left(x_i\xch_i+\sum\nolimits_{j<i} s_{ji}\right)^d
\in\+_{k\leq d}\X_{(k)}\cdot \F\W\cdot \Xch_{(k)}.$$
Since $M$ is locally nilpotent for $\Xch_+$, there exists
$K\in\Z_{\geq 0}$ such that $\Xch_{(k)}v=0$ if $k>K$.
Hence for any $i\in[1,n]$ and $d\in\Z_{\geq 0}$, it holds that
$\ech_i^d v$ is contained in
$\+_{k\leq K}\X_{(k)}\cdot \F\W\cdot \Xch_{(k)}v$,
which is finite-dimensional.
This implies $\dim(\UU v)<\infty$.

By Corollary~\ref{cor;pi_and_H},
any element in $\daffH\*_{\ratH}M$ is of the form
$\pi^k \bar v$ for some $k\in\Z_{\leq 0}$ and $v\in M$,
where $\bar v$ denotes the image of $1\*v$ in $\daffH\*_{\ratH}M$.
Since $\pi \UU\pi^{-1}=\UU$,
we have
$\dim\left( \UU \pi^k\bar v\right)
=\dim ( \UU v)<\infty.$
\qed

\medskip
Let $\Res(-)$ denote the restriction functor
$\daffHmod\to\ratHmod$.
\begin{lemma}\label{lem;inclusion}
$\rm{(i)}$
Let $M$ be an object of $\O(\ratH)$.
Then the natural embedding
$M\hookrightarrow \daffH\*M$
induces an injective $\ratH$-homomorphism
$$
M\hookrightarrow
\Res
\left(\In M
\right).$$

\smallskip\noindent
$\rm{(ii)}$
Let $N$ be an object of $\O(\daffH)$.
Then the product map
$\daffH\*N\to N$ induces
an $\daffH$-isomorphism
$$\daffH\otimes_{\ratH}({\Res N})\overset{\sim}{\to} N.$$
\end{lemma}
\noindent{\it Proof.}
(i) It is obvious that the map is an $\ratH$-homomorphism. 
Let us prove that it is injective.
Let $v\in M\setminus\{0\}$.
Put $\omega=\sum_{i\in[1,n]}\xch_i\in\Xch_{(1)}\subset\ratH$.
Then
we have $[\omega^p, x_i]= p\kappa\omega^{p-1}$
for any $p\in\Z_{\geq1}$.
Since $M$ is locally nilpotent for $\Xch_+$,
there exists $p_0\in\Z_{\geq 1}$
such that $\omega^{p_0-1}v\neq 0$
and $\omega^{p_0}v=0$.
For any $i\in[1,n]$,
we have $\omega^{p_0} x_iv={p_0}\kappa\omega^{p_0-1}v\neq 0$,
and hence $x_iv\neq 0$.
This implies $fv\neq 0$ for any monomial $f$ in $x_1,\dots,x_n$.
This proves the statement. 

\noindent
(ii) Follows easily from the isomorphism
 $\daffH\otimes_{\ratH}N\cong \F[\pi^{\pm1}]\*_{\F[\pi]}N$.
\qed
%
\section{The adjoint functor $\nil$}
Now we introduce a functor $\O(\daffH)\to\O(\ratH)$,
which turns out to be the right adjoint functor of $\induction.$

For an $\daffH$-module $M$,
denote by $\nil(M)$ the subspace of $M$ consisting of
the $\Xch_+$-nilpotent elements.
\begin{lemma}
 If $M$ is in $\O(\daffH)$,
then $\nil(M)$ is an $\ratH$-submodule
of $\Res M$ and it is finitely generated over $\ratH$.
\end{lemma}
\noindent{\it Proof.}
Let $v\in \nil(M)$. Then
it is obvious that $s_iv$ and $\xch_i v$ are in
$\nil(M)$.
It is easily shown by induction that
$\Xch_{(p)} x_i\subseteq \X_{(1)}\Xch_{(p)}+\F\W  \Xch_{(p-1)}$
for any $p\in\Z_{\geq 1}$.
This implies $x_iv\in  \nil(M)$
and hence $\nil(M)$ is  an $\ratH$-submodule.

Let us prove that $\nil(M)$ is finitely generated.
Take a sequence
\begin{equation}
  \label{eq;nseq}
  0=N^{(0)}\subsetneq N^{(1)}\subsetneq
\dots\subsetneq N^{(k)}\subsetneq \cdots
\end{equation}
of finitely generated $\ratH$-submodules of $\nil(M)$.
Since each $N^{(k)}$ is in $\O(\ratH)$,
it follows from Lemma~\ref{lem;inclusion}-(i)
that $\daffH\*_{\ratH}N^{(k)}\subsetneq
\daffH\*_{\ratH}N^{(k+1)}$ for all $k$.
By the exactness of the induction functor on $\ratHmod$ and by
 Lemma~\ref{lem;inclusion}-(ii),
we have
 $\daffH\*_{\ratH}\nil(M)\subseteq\daffH\*_{\ratH}\Res(M)
\cong M$.
Hence  $\daffH\*_{\ratH}\nil(M)$
is in $\O(\daffH)$
and its composition series is finite. 
Therefore the sequence \eqref{eq;nseq} must end in finite steps,
and we have $N^{(K)}=\nil(M)$ for some $K\in\Z_{\geq 0}$.
\qed

\medskip
Therefore we obtain the functor
$\nil:\O(\daffH)\to\O(\ratH)$.
\begin{corollary}\label{cor;rightadjoint}
The functor $\nil:\O(\daffH)\to\O(\ratH)$ is the right adjoint functor 
of the functor $\induction: \O(\ratH)\to\O(\daffH);$ namely, 
we have a natural isomorphism
 \begin{equation}
   \Hom_{\daffH}(\In M,N)
\cong\Hom_{\ratH}(M,\nil(N))
 \end{equation}
  for any $M$ in $\O(\ratH)$ and any $N$ in $\O(\daffH)$.
\end{corollary}
\begin{proposition}\label{pr;nilind=id}
  Let $M$ be in $\O(\ratH)$.
Then there exists a
natural $\daffH$-isomorphism $\nil(\In M)\cong M$.
\end{proposition}
\noindent{\it Proof.}
Using Lemma~\ref{lem;inclusion}-(i),
we have
the following 
sequence of  embeddings of $\ratH$-modules:
\begin{equation}\label{eq;ratemb}
M\overset{i}{\hookrightarrow} \nil(\In M)\hookrightarrow
\Res(\In M).
\end{equation}
By applying $\induction$, we have embeddings of
$\daffH$-modules
\begin{equation}\label{eq;daffemb}
\widetilde M\overset{\tilde i}{\hookrightarrow}
 \In{(\nil(\widetilde M))}\hookrightarrow\In{(\Res \widetilde M)},
\end{equation}
where we put $\widetilde M=\In M$.
By combining it with
the isomorphism
$\In{(\Res \widetilde M)}\overset{\sim}{\to} \widetilde M$
given in Lemma~\ref{lem;inclusion}-(ii),
we have
an $\daffH$-homomorphism
$\widetilde M\to \widetilde M$,
and it easily follows that this map equals the identity map.
Therefore the embedding $\tilde i$ in \eqref{eq;daffemb}
is an isomorphism, and hence
 we have  
$\In{(\Coker\, i)}=\Coker\, \tilde i=0.$
This implies $\Coker\,i=0$ 
by Lemma~\ref{lem;inclusion}-(i).
Namely, the embedding $i:M\to \nil(\In M)$ in \eqref{eq;ratemb}
is an isomorphism.
\qed
\begin{theorem}\label{th;fullyfaithful}
 The functor
$\induction:\O(\ratH)\to\O(\daffH)$
is exact and fully faithful.
\end{theorem}
\noindent
{\it Proof.}
The exactness has been shown in 
Corollary~\ref{cor;flat}.
The claim that $\induction$
is fully faithful
follows from
 Corollary~\ref{cor;rightadjoint}
and Proposition~\ref{pr;nilind=id}
\qed
\begin{remark}
It is conjectured 
in \cite[Remark 5.17]{GGOR}
that the category $\O(\ratH)$
is equivalent to the category of finite dimensional
modules over the $q$-Schur algebra with  
$q={\rm exp}(\frac{2\pi i}{\kappa}).$
The analogous statement for $\daffH$
has been proved in \cite[Theorem 5.2]{VV}, where they
gave a categorical equivalence between
$\O(\daffH)$ and 
the category of finite dimensional modules
over the affine $q$-Schur algebra.
\end{remark}
\section{Standard and irreducible modules}
Introduce the subalgebra
 $\ratHp=\F\W\cdot \Xch$ of $\ratH$.
Let $\phi:\ratHp\to\F\W$ be
the algebra homomorphism
given by
\begin{equation*}
 \phi(s_i)=s_i\ (i\in[1,n-1]),\quad
\phi(x_i)=0\ (i\in[1,n]).
\end{equation*}
For  a finite-dimensional $\W$-module $E$,
let ${}^\phi\! E$ denote the $\ratHp$-modules
given by the composition
$\ratHp\stackrel{\phi}{\to}\F\W\to\End(E)$.

We define the standard module of $\ratH$ associated with $E$ by
\begin{equation*}
\ratst(E)=\ratH\otimes_\ratHp \!{}^\phi\! E.
\end{equation*}
There exists an algebra homomorphism
$\ev:\affH\to \F\W$ given by
\begin{equation*}
\ev(s_i)= s_i\ (i\in[1,n-1]),\quad 
\ev(\ech_i)= \sum\nolimits_{j<i}s_{ji}\ (i\in[1,n]),
\end{equation*}
and define the $\affH$-module
${}^\ev\! E$ similarly as ${}^\phi\!E$.
We define the standard module of $\daffH$ associated with $E$ by
\begin{equation*}
  \daffst(E)=\daffH\otimes_\affH\! {}^\ev\! E.
\end{equation*}
Observe that the standard modules $\ratst(E)$ (resp. $\daffst(E)$)
are objects of $\O(\ratH)$ (resp. $\O(\daffH)$).
\begin{theorem}\label{th;ind_O}
${\rm (i)}$  Let $E$ be a finite-dimensional $\W$-module.
Then $$\daffH\otimes_\ratH \ratst(E)\cong \daffst(E).$$

\smallskip
\noindent
${\rm (ii)}$ Let $M$ be in $\O(\ratH)$.
Then  $\daffH\*_\ratH M$ is irreducible if and only if
$M$ is irreducible.
\end{theorem}
\noindent{\it Proof.}
(i) Observe that
 $\xch_i$ acts trivially on
the image of  $\!{}^\ev\! E$ in $\daffst(E)$.
Thus
we have
an $\ratH$-homomorphism
$\ratst(E) \to\daffst(E)$,
and hence
a surjective $\daffH$-homomorphism
$\daffH\*_{\ratH}\ratst(E)\to\daffst(E)$.
This is isomorphism since both sides are isomorphic
to $\F P\otimes E$ as $\F P$-modules.

\noindent
(ii) The statement follows immediately from
Theorem~\ref{th;fullyfaithful}.
\qed

\medskip
Let $E$ be a finite-dimensional irreducible $\W$-module.
A similar argument as in
the theory of highest weight modules 
for Lie algebras 
implies that the standard module $\ratst(E)$
has a unique simple quotient,
and any simple object in $\O(\ratH)$ is
given as such a quotient \cite[Proposition 2.11]{GGOR}.
\begin{corollary}$($cf. \cite[Proposition 2.5.3]{AST}
\cite[Theorem 8.2]{Su;classification}$)$
  Let $E$ be a finite-dimensional  irreducible
$\W$-module.
Then $\daffst(E)$ has a unique simple quotient module.
\end{corollary}
We denote by $\ratL(E)$ and  $\daffL(E)$
the unique simple quotient of $\ratst(E)$
and  $\daffst(E)$ respectively.
As a direct consequence of Theorem~\ref{th;ind_O},
we have $\daffH\*_{\ratH}\ratL(E)\cong
\daffL(E)$,
\begin{corollary}\label{cor;multiplicity}
 Let $E$ and $F$ be finite-dimensional  irreducible
$\W$-modules. Then
$$ [\ratst(E):\ratL(F)]=[\daffst(E):\daffL(F)],$$
where $[M:N]$ denotes the multiplicity of the irreducible module
$N$ in the composition series of $M$.
\end{corollary}
\begin{remark}
The category $\O(\daffH)$ has been also related to
the category $\O$ for the (original) double affine Hecke algebra 
(\cite[Proposition 2.5]{VV}),
and, as a consequence, the multiplicity for $\daffH$
in Corollary~\ref{cor;multiplicity} equals
the corresponding multiplicity for the double affine Hecke algebra.
The latter can be expressed in terms of
the Kazhdan-Lusztig polynomials (\cite[Theorem 8.5]{Va}).
\end{remark}
\section{$\UU$-semisimple  modules}
Recall that $\ratH$ includes the degenerate affine Hecke algebra
$\affH$ 
and in particular the commutative subalgebra $\UU$.

For an $\affH$-module $M$
and $\zeta=(\zeta_1,\dots.\zeta_n)\in \wt$,
set
$$M_\zeta=\{v\in M\mid \ech_iv=\zeta_i v\ 
\hbox{ for all }i\in[1,n]\}.$$
We call an element of $M_\zeta$ a weight vector of weight $\zeta$.

Define the category $\Oss(\daffH)$
as the full subcategory of
$\O(\daffH)$
consisting of 
$\UU$-semisimple modules, i.e., an object 
$M$ in $\O(\daffH)$ is in $\Oss(\daffH)$
if $M=\oplus_{\zeta\in\wt}M_\zeta$.
Similarly, define the category $\Oss(\ratH)$
as the full subcategory
 of $\O(\ratH)$
consisting of $\UU$-semisimple modules.

It follows that $\dim M_\zeta<\infty$ for any $\zeta\in\F^n$
and for any $M\in \Oss(\daffH)$ or $M\in \Oss(\ratH)$.
%
\begin{proposition}\label{pr;OsstoOss}
An $\ratH$-module $M$ is in $\Oss(\ratH)$
if and only if $
\In M$
is in $\Oss(\daffH)$.
\end{proposition}
%
\noindent{\it Proof.}
Suppose that  $\In M$ is in $\Oss(\daffH)$.
Then
$M\cong\nil(\In M)\subseteq \In M$,
and it is in $\Oss(\ratH)$.

Suppose that  $M$ is in $\Oss(\ratH)$.
Since
$\daffH\*_{\ratH}M=\F[\pi,\pi^{-1}]\*_{\F[\pi]} M$
by Corollary~\ref{cor;pi_and_H},
any element in $\daffH\*_{\ratH}M$ is of the form
$\pi^k  v$ for some $k\in\Z$ and $v\in M$.
The relation \eqref{eq;pi_ech} implies that
if $v\in M$ is a weight vector for $\ech_1,\dots,\ech_n$
then $\pi^k v$ is also a weight vector.
This implies that $\In M$ is in $\Oss(\daffH)$.
\qed

\medskip
Let $\m\in\Z_{\geq1}$.
Let $\part(\m,n)$ denote the set of partitions of $n$ with $\m$ 
nonzero components:
\begin{align*}
&\part(\m,n)=\\
&\left\{\lm=(\lm_1,\lm_2,\dots,\lm_\m)\in\Z^\m\
\left\vert\,\lm_1\geq\lm_2\geq\dots\geq\lm_\m>0,\
\sum_{i\in[1,\m]}\lm_i=n\right.\right\}.
\end{align*}
For $\lm\in\part(\m,n)$,
let $E_\lm$ denotes the irreducible $\W$-module
corresponding to the partition
(or the Young diagram) $\lm$.
We write $\ratL(\lm)$ and $\daffL(\lm)$ for
$\ratL(E_\lm)$ and $\daffL(E_\lm)$ respectively.

Let $\Irr(\Oss(\ratH))$ (resp. $\Irr(\O(\daffH)$) denote
the set of isomorphism classes
of the irreducible modules in $\Oss(\ratH)$
(resp. $\O(\daffH)$).
Recall that the  assignment $\lm\mapsto \ratL(\lm)$ gives
a one-to-one correspondence
$\bigsqcup_{\m\in[1,n]}\part(\m,n)\isomto\Irr(\O(\ratH))$
(\cite[Section 2.5.1]{GGOR}).

Now let $\kappa\in\Z$, and set
\begin{equation*}
\part_\kappa(\m,n)=
\{\lm\in \part(\m,n)\mid
\kappa-\m-\lm_1+\lm_\m\in\Z_{\geq0}\}.
\end{equation*}
Proposition~\ref{pr;OsstoOss} and
the classification of the irreducible modules in $\Oss(\daffH)$
(\cite[Theorem 6.5]{Ch;fourier} \cite[Corollary 4.23]{SV})  implies
the following:
\begin{theorem}\label{th;classification}
${\rm (i)}$ For $\kappa\in\Z_{>0}$,
the assignment $\lm\mapsto \ratL(\lm)$ gives
a one-to-one correspondence
$$\bigsqcup_{\m\in[1,n]}\part_\kappa(\m,n)\isomto\Irr(\Oss(\ratH)).$$

\smallskip
\noindent
${\rm (ii)}$ For $\kappa\in\Z_{<0}$,
the assignment $\lm\mapsto \op{\mathcal L}^\rat_{-\kappa}(\lm)$ gives
a one-to-one correspondence
$$\bigsqcup_{\m\in[1,n]}\part_{-\kappa}(\m,n)
\isomto\Irr(\Oss(\ratH)),$$
where
$\op{\mathcal L}^\rat_{-\kappa}(\lm)$ denotes the
$\ratH$-module given by the composition
$\ratH
\stackrel{\sigma}{\to} \H_{-\kappa}
\to \End({\mathcal L}^\rat_{-\kappa}(\lm))$.
\end{theorem}
\begin{remark}
We stated
Theorem~\ref{th;classification} 
for $\kappa\in\Z$
since we treated mainly this case in \cite{SV},
but most of the results obtained there can be generalized
for general $\kappa\in\F\setminus\{0\}$.
In particular it follows that
Theorem~\ref{th;classification} can be generalized by
setting $\part_\kappa(\m,n)=\part(\m,n)$ for 
$\kappa\in\F\setminus\Q$, and
\begin{equation*}
\part_\kappa(\m,n)=
\{\lm\in \part(\m,n)\mid
r-\m-\lm_1+\lm_\m\in\Z_{\geq0}\}
\end{equation*}
for $\kappa=r/p\in\Q_{>0}$ with $r,p\in\Z_{>0}$, $(r,p)=1$
(cf. \cite[Theorem 6.5]{Ch;fourier}).
\end{remark}
\begin{remark}
  For $\lm\in\part(\m,n)$, it holds that
$\op{\mathcal L}^\rat_{-\kappa}(\lm)\cong \ratL(\lm')$,
where $\lm'$ is the conjugate of $\lm$:
$$\lm'_i=\sharp\{a\in[1,\m]\mid \lm_a\geq i\}\quad
(i\in[1,\lm_1]).$$
\end{remark}
\begin{remark}
It can be seen that
{\it the Knizhnik-Zamolodchikov functor}
investigated in \cite{GGOR}
(and in \cite{VV} for $\daffH$)
transforms the irreducible representations
$\ratL(\lm)$ for $\lm\in\part_\kappa(\m,n)$ to Wenzl's
representations~\cite{Wenzl}
of the Hecke algebra. 
\end{remark}

\end{document}